\documentclass[journal]{IEEEtran}

\usepackage{cite}
\usepackage{amsmath,amssymb,amsfonts}
\usepackage{color}
\usepackage{algorithmic}
\usepackage{graphicx}
\usepackage{textcomp}
\usepackage[english]{babel}
\usepackage[utf8]{inputenc}
\usepackage{multirow}
\usepackage{tabularx}
\usepackage{nomencl}
\usepackage{kantlipsum}
\allowdisplaybreaks
\makenomenclature
\def\BibTeX{{\rm B\kern-.05em{\sc i\kern-.025em b}\kern-.08em
    T\kern-.1667em\lower.7ex\hbox{E}\kern-.125emX}}
\usepackage{etoolbox}
\renewcommand\nomgroup[1]{%
  \item[\bfseries
  \ifstrequal{#1}{V}{Variables}{%
  \ifstrequal{#1}{P}{Parameters}{%
  \ifstrequal{#1}{S}{Sets}{}}}%
]}

\ifCLASSINFOpdf

\else

\fi

\begin{document}

\title{Micro Water-Energy Nexus: Optimal Demand-Side Management and Quasi-Convex Hull Relaxation}

\author{Qifeng~Li,~\IEEEmembership{Member,~IEEE,}
        Suhyoun~Yu, Ameena S. Al-Sumaiti, 
        and~Konstantin~Turitsyn,~\IEEEmembership{Member,~IEEE}
\thanks{This work is supported by the MI-MIT Cooperative Program under grant MM2017-000002.}
\thanks{Q. Li, S. Yu, and K. Turitsyn are with the Department
of Mechanical Engineering, Massachusetts Institute of Technology, Cambridge,
MA, 02139 USA, e-mail: \{qifengli, syu2, turitsyn\}@mit.edu.}
\thanks{A. Al-Sumaiti is with the Department
of Electrical Engineering \& Computer Science, Masdar Institute, Khalifa University of Science and Technology, Abu Dhabi, UAE, e-mail: aalsumaiti@masdar.ac.ae}}

\markboth{May~2018}%
{Shell \MakeLowercase{\textit{et al.}}: Bare Demo of IEEEtran.cls for IEEE Journals}

\maketitle

\begin{abstract}
This paper investigates the water network's potential ability to provide demand response services to the power grid under the framework of a distribution-level water-energy nexus (micro-WEN). In particular, the hidden controllability of water loads, such as irrigation systems, was closely studied to improve the flexibility of electrical grids. A optimization model is developed for the demand-side management (DSM) of micro-WEN, and the simulation results assert that grid flexibility indeed benefits from controllable water loads. Although the proposed optimal DSM model is an intractable mixed-integer nonlinear programming (MINLP) problem, quasi-convex hull techniques were developed to relax the MINLP into a mixed-integer convex programming (MICP) problem. The numerical study shows that the quasi-convex hull relaxation is tight and that the resulting MICP problem is computationally efficient.

\end{abstract}


\IEEEpeerreviewmaketitle

\nomenclature[V]{$P_{ij,t}$}{Active branch power in line $ij$}
\nomenclature[V]{$Q_{ij,t}$}{Reactive branch power in line $ij$}
\nomenclature[V]{$\mathcal{V}_{i,t}$}{Square of voltage at bus $i$}
\nomenclature[V]{$\mathcal{I}_{ij,t}$}{Square of current in line $ij$}
\nomenclature[V]{$P_{i,t}^G $}{Active power of diesel generator at bus $i$}
\nomenclature[V]{$Q_{i,t}^G $}{Reactive power of diesel generator at bus $i$}
\nomenclature[V]{$P_{i,t}^{ES} $}{Active power of the BESS at bus $i$}
\nomenclature[V]{$ Q_{i,t}^{ES}$}{Reactive power of the BESS at bus $i$}
\nomenclature[V]{$P_{i,t}^{\mathrm{Pump}} $}{ Active consumptions of the pump at bus $i$}
\nomenclature[V]{$ Q_{i,t}^{\mathrm{Pump}}$}{Reactive consumptions of the pump at bus $i$}
\nomenclature[V]{$L_{i,t}^{ES} $}{ Active power loss of the BESS at bus $i$}
\nomenclature[V]{$x_t $}{ Vector of water flow in pipes}
\nomenclature[V]{$y_t $}{ Vector of water head at nodes}
\nomenclature[V]{$w_{i,t}^{UT} $}{ Water flow to the utility-owned water tank at node $i$}
\nomenclature[V]{$w_{i,t}^{CT} $}{ Water flow to the customer-owned water tank at node $i$}
\nomenclature[V]{$w_{i,t}^G $}{ Water flow injected from the water source in pipe $i$}
\nomenclature[V]{$y_{k,t}^G $}{ Head gain imposed by the pump in pipe $k$}
\nomenclature[V]{$\alpha_{i,t}$}{Binary variable denoting the on/off state of the pump in pipe $i$}
\nomenclature[P]{$A $}{Incidence matrix of water network}
\nomenclature[P]{$B_{i}$}{Coefficient of pump characteristics in pipe $i$}
\nomenclature[P]{$C$}{Cost function of the whole system}
\nomenclature[P]{$P_{i,t}^L $}{Active electric load at bus $i$}
\nomenclature[P]{$Q_{i,t}^L$}{Reactive electric load at bus $i$}
\nomenclature[P]{$P_{i,t}^{RE} $}{Active power output of renewable at bus $i$}
\nomenclature[P]{$Q_{i,t}^{RE}$}{Reactive power output of renewable at bus $i$}
\nomenclature[P]{$r_{ij}$}{ Resistance of line $ij$}
\nomenclature[P]{$x_{ij}$}{ Reactance of line $ij$}
\nomenclature[P]{$r_i^{\mathrm{Batt}}$}{ Loss coefficient related to the battery of BESS}
\nomenclature[P]{$r_i^{Cvt}$}{ Loss coefficient related to the converter of BESS}
\nomenclature[P]{$E_{i,0}^{ES}$}{ Initial state of charging of BESS at bus $i$}
\nomenclature[P]{$h_i $}{Elevation at node $i$ of the water network}
\nomenclature[P]{$R_{k}^P$}{ Head loss coefficient of pipe $k$}
\nomenclature[P]{$c_t$}{Locational marginal price at PCC}
\nomenclature[P]{$c_{1,i}$}{Coefficient of the first-order term in the cost function of diesel generator at bus $i$}
\nomenclature[P]{$c_{2,i}$}{Coefficient of the second-order term in the cost function of diesel generator at bus $i$}
\nomenclature[S]{$\mathcal{E}_E$}{Edge set of the electricity network}
\nomenclature[S]{$ \mathcal{E}_W$}{Edge set of the water network}
\nomenclature[S]{$\mathcal{N}_E$}{Bus set of the electricity network}
\nomenclature[S]{$ \mathcal{N}_W$}{Node set of the water network}
\nomenclature[S]{$\mathcal{E}_W^P$}{Set of pipes with a pump installed}
\nomenclature[S]{$\mathcal{N}_E^S$}{Set of bus with a BESS connected}
\nomenclature[S]{$ \mathcal{N}_W^S$}{Set of nodes connected to a tank}
\nomenclature[S]{$\mathcal{N}_E^G$}{Set of bus with controllable generations}
\nomenclature[S]{$\mathcal{N}_E^P$}{Set of bus with a pump connected}

\printnomenclature

\section{Introduction}

\IEEEPARstart{M}{odern} day water and power systems are closely intertwined as a coupled system, commonly referred to as the water-energy nexus (WEN) \cite{LiModeling, SiddiqiWater, SanthoshRealtime, ZhangEnergy, SanthoshOptimal}. On one hand, most of the water facilities consume electrical energy. For instance, ground water pumping and seawater desalination account for roughly 12$\%$ of the total electric energy consumption in the Arabian Gulf regions \cite{SiddiqiWater}. On the other hand, water usage is inevitable in refining fuels and generating electric power. Despite the two networks' inevitable interdependency, water and energy networks have traditionally been operated independently from one another, and the idea of co-operating the two in parallel has long been glossed over.



In recent years, researchers have started to direct their attention to water system's potential ability to provide demand response (DR) services \cite{MenkeDemonstrating}. It is believed that higher cost-efficiency can be achieved by co-operating the water and power systems\footnote{California water utilities stake out new role in energy programs to finance their future, available at http://artemiswaterstrategy.com/slopgepartnership/.}. In 2016, PG$\&$E built up an efficiency partnership with the water utility in the city of San Luis Obispo which is the first-of-its-kind\footnote{First-of-its-kind efficiency partnership with PG$\&$E expected to save SLO millions in energy bills, available at http://kcbx.org/post/first-its-kind-efficiency-partnership-pge-expected-save-slo-millions-energy-bills$\#$stream/0.}. The power sector suffers from an unprecedented amount of network uncertainty with the ever-increasing penetration of intermittent renewable energy and electro-mobility \cite{PalenskyDemand}, and the idea of exploiting the flexibility of water network has come under the spotlight as a possible solution. The co-operation of two systems allows the water system to fast and accurately response to the energy imbalance caused by the uncertain renewable energy generation or even contingencies on the electricity side. With this solution, the overall security and reliability of water and power systems will benefit from the nexus operation mode.

\begin{figure*}[th!]
\centering
\includegraphics[width=0.8\textwidth]{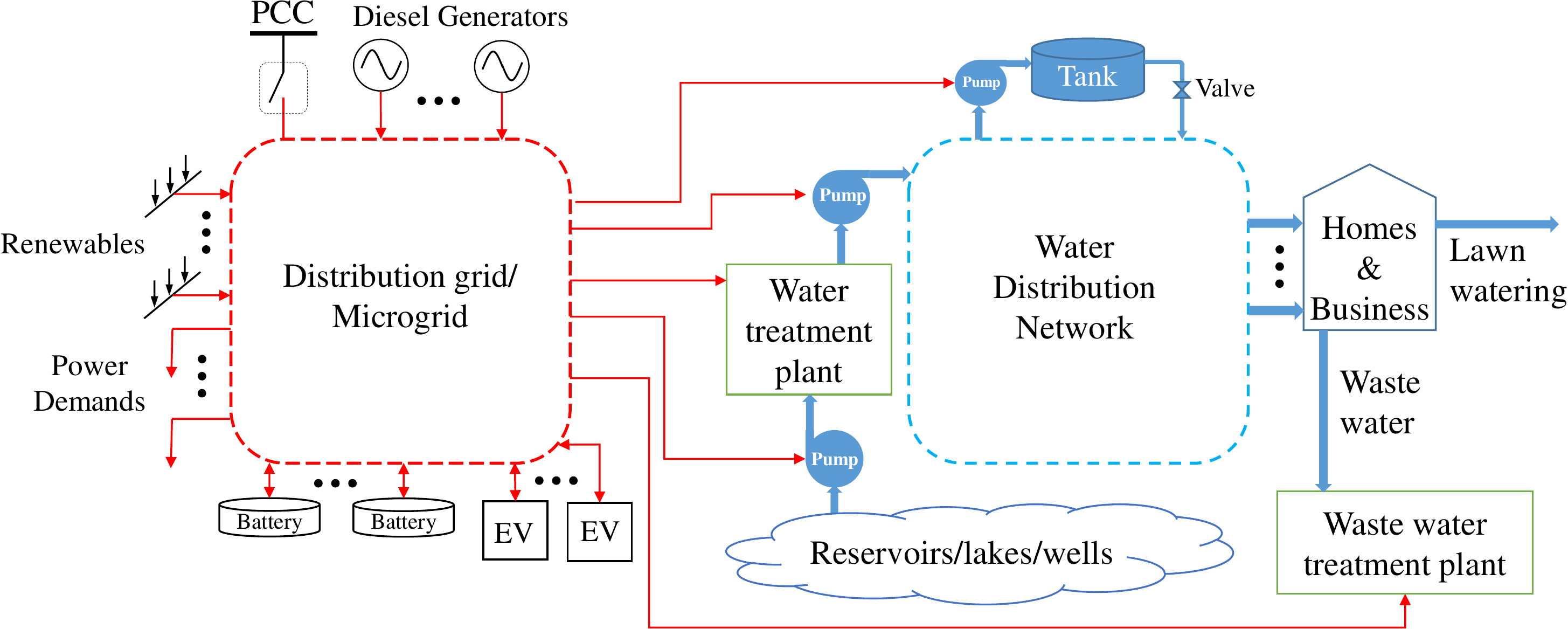}
  \caption{A physical structure of the micro water-energy nexus.}
  \label{fig:WEN}
\end{figure*}


This paper focuses on developing a mathematical tool to assess the flexibility and responsiveness of a given micro-WEN, wherein both the water and energy networks are at distribution levels, to the energy imbalance. The AC power flow integrated with battery energy storage systems (BESSs) and high penetration of renewable energy sources \cite{LiConvex} is used in the developed  mathematical model. The water pipe network model is characterized by the directed Darcy-Weishach equation \cite{ChaudhryApplied} and allow for water flow directions to reverse, and the on/off status of pumps are represented by integer variables. The resulting mathematical model for the micro-WEN is unfortunately a nonlinear mixed-integer model. To mitigate its intractability, a quasi-convex hull relaxation \cite{LiThe,LiConvex} that is sufficiently tight and efficient is developed for the micro-WEN model.

Much of the micro-WEN's flexibility will rise from the water sector, as electricity-driven water services---pumping, cooling, desalination, water treatment---are time-flexible. Further flexibility can be introduced by including controllable water loads such as irrigation services, which this paper uses to demonstrate hidden DR capabilities of the water network. The flexibility of pumps is constrained by water tank capacities that are physically limiting by nature, but incorporating controllable discharge scheme for water tanks should resolve the issue. Moreover, pumps, tanks and irrigation systems are to be jointly used to create virtual energy storages to power grids to alleviate the stress on the WEN in case of limited energy supply.

\section{Problem Formulations}
\subsection{Physical Model of Micro-WEN}

The schematic of the micro-WEN is given in Figure \ref{fig:WEN}. The electricity side is a distribution network, or a micro-grid, integrated with renewable energy and BESSs. The water side consists of a pipe network, pumps, utility- and customer-owned tanks and irrigation systems. EV's and water treatment facilities---including desalination, water and waste water treatment and recycling---are not considered for the sake of simplicity, but those two elements can be easily incorporated and will be considered in future research.


It has been discussed in \cite{LiModeling} that such a micro-WEN is a fundamental infrastructure of a smart building/city/village. Under the environment of smart buildings\cite{Morvaj}/cities\cite{Albino}/villages\cite{Heap}, all infrastructures will be connected through the emerging Internet of Things techniques. In order to operate and control such a connected physical system, it is essential to develop the mathematical model and design the optimization algorithms for optimal resource allocation. 

\subsection{Mathematical Model of Micro-WEN}
This section introduces a multi-period mathematical model of the micro-WEN. Unless otherwise stated, the subscript $t$ denotes the discrete time period. The structure of an AC-microgrid or an electric distribution system is usually radial. Consequently, we use a \textit{DistFlow} model \cite{BaranOptimal} integrated with renewable generation and batteries  to describe the electricity network. The detailed model is given as
%
\begin{subequations} \label{PF}
\begin{gather} 
\begin{split}
&P_{i,t}^\mathrm{G}+P_{i,t}^\mathrm{RE}-P_{i,t}^{\mathrm{Pump}}-P_{i,t}^\mathrm{L}+P_{i,t}^\mathrm{ES} \\
&= r_{ij} \mathcal{I}_{ji,t}-P_{ji,t}+\sum_{k \in \mathcal{D}_i}P_{ik,t} 
\end{split}\\
\begin{split}
&Q_{i,t}^\mathrm{G}+Q_{i,t}^\mathrm{RE}-Q_{i,t}^{\mathrm{Pump}}-Q_{i,t}^L+Q_{i,t}^\mathrm{ES} \\
&= x_{ij} \mathcal{I}_{ji,t}-Q_{ji,t}+\sum_{k \in \mathcal{D}_i}Q_{ik,t}
\end{split}\\
\mathcal{V}_{i,t}-\mathcal{V}_{k,t}+(r_{ik}^2+x_{ik}^2)\mathcal{I}_{ik,t}=2(r_{ik}P_{ik,t}+x_{ik}Q_{ik,t})\\
P_{ik,t}^2+Q_{ik,t}^2 = \mathcal{V}_{i,t}\mathcal{I}_{ik,t}, \  (ik \in \mathcal{E}_\mathrm{E}) \\
P_{ik,t}^2+Q_{ik,t}^2 \le \overline{S}_{ik}^2, \  (ik \in \mathcal{E}_\mathrm{E}) \\
0 \le \mathcal{I}_{ik,t} \le \overline{\mathcal{I}}_{ik}, \\
\underline{\mathcal{V}}_i \le \mathcal{V}_{i,t} \le \overline{\mathcal{V}}_i \\
 \underline{P}_i^\mathrm{G},\underline{Q}_i^\mathrm{G} \le P_{i,t}^\mathrm{G},Q_{i,t}^\mathrm{G} \le \overline{P}_i^\mathrm{G},\overline{Q}_i^\mathrm{G},
\end{gather}
\end{subequations}
where $i \in \mathcal{N}_E$ and $ik \in \mathcal{E}_E$. For the sake of simplicity, this paper only consider the balanced cases. It has been proved in literature that the convex relaxations of the \textit{DistFlow} for balanced networks can be easily extended to the unbalanced cases under some mild approximations. Therefore, the proposed approach in this paper can be easily leveraged to the cases of three-phase unbalanced distribution networks or microgrids.

The following nonlinear model of a battery energy storage unit is incorporated into the overall mathematical model of the micro-WEN. Please refer to \cite{LiConvex} for more details about this BESS model. For $\forall i \in \mathcal{N}_E^S$, we have
\begin{subequations} \label{BESS}
\begin{gather} 
(r_i^{\mathrm{Batt}} + r_i^{\mathrm{Cvt}})(P_{i,t}^\mathrm{ES})^2 + r_i^\mathrm{Cvt}(Q_{i,t}^\mathrm{ES})^2 = L_{i,t}^\mathrm{ES}\mathcal{V}_{i,t} \\
(P_{i,t}^\mathrm{ES})^2 + (Q_{i,t}^\mathrm{ES})^2 \leqslant (\overline{S}_{i}^\mathrm{ES})^2 \\
\underline{E}_{i}^\mathrm{ES} \leqslant E_{i,0}^\mathrm{ES} - \sum^t_{t=0} (P_{i,t}^\mathrm{ES}+L_{i,t}^\mathrm{ES}) \leqslant \overline{E}_{i}^\mathrm{ES}.
\end{gather}
\end{subequations}

We make the following assumptions for the water distribution networks: (1) The pipe network is a directed graph $\mathcal{G}_\mathrm{W}=(\mathcal{N}_\mathrm{W},\mathcal{E}_\mathrm{W})$ with incidence matrix $A$ such that  $A_{ik} \in \{-1,\ 0,\ 1\}$ for all $i,k$; (2) A pump is considered as a type of pipe that imposes a head gain when the pump is on and closed otherwise; (3) The pump converts the electric power into a mechanical power at a constant efficiency of $\eta$; (4) The power factors of pumps are fixed, namely $P_{k,t}^{\mathrm{Pump}}/Q_{k,t}^{\mathrm{Pump}}$ is constant. The resulting model can be expressed as:
\begin{subequations} \label{WDS}
\begin{gather} 
\sum_{k \in \mathcal{E}_W}A_{ik}f_{k,t} = f_{i,t}^\mathrm{G} - f_{i,t}^\mathrm{UT} - f_{i,t}^\mathrm{CT},(i \in \mathcal{N}_\mathrm{W}) \\
y_{i,t}-y_{j,t}+h_i-h_j= R_{k}^\mathrm{P}\mathrm{sgn}(f_{k,t})f_{k,t}^2,(k \in \mathcal{E}_\mathrm{W}\setminus \mathcal{E}_\mathrm{W}^\mathrm{P})\\
\begin{cases}
\begin{array}{ll}
 \begin{split}
 &y_{i,t}-y_{j,t}+h_i-h_j\\&+y_{k,t}^\mathrm{G}= R_k^\mathrm{P}f_{k,t}^2
 \end{split} & \text{if} \ \alpha_{k,t}=1 \\
 f_{k,t} = 0, & \text{if} \ \alpha_{k,t}=0
\end{array}, (k \in \mathcal{E}_\mathrm{W}^\mathrm{P}) 
\end{cases} \\
y_{k,t}^\mathrm{G} = B_{k}f_{k,t}+C_{k} \quad (k \in \mathcal{E}_\mathrm{W}^\mathrm{P}) \label{eq:yf}\\
 \underline{S}_{i}^\mathrm{w} \leqslant S_{i,0}^\mathrm{w} + \sum^t_{t=0} w_{i,\tau}^\mathrm{UT} \leqslant \bar{S}_{i}^\mathrm{w},\quad (i \in \mathcal{N}_\mathrm{w}^\mathrm{S})\\
 \underline{f} \le f_t \le \overline{f}, \\
 \underline{y} \le y_t \le \overline{y}, \\
 \underline{w}_i^\mathrm{G} \le w_{i,t}^G \le \overline{w}_i^\mathrm{G}, \  (i \in \mathcal{N}_W^\mathrm{G}) \\
 \underline{w}_i^\mathrm{S} \le w_{i,t}^\mathrm{S} \le \overline{w}_i^\mathrm{S}, \  (i \in \mathcal{N}_\mathrm{W}^\mathrm{S})
\end{gather}
\end{subequations}
where $A$ is a $|\mathcal{N}_\mathrm{W}|\times|\mathcal{E}_\mathrm{W}|$ incidence matrix and pipe $k$ connects nodes $i$ and $j$. Equation (3a) represents the mass balance of the water network; constraints (3b) and (3c) formulate the hydraulic characteristics of a normal pipe and the pipe with a pump installed respectively; constraint (3e) denotes the state of charging of the water tanks; (3f)-(3i) are system constraints; $\mathrm{sgn}(f)=-1$ if $f \le 0$ or, otherwise $\mathrm{sgn}(f)=1$. When $\alpha_{k,t}=1$, the quantity $f_{k,t}$ in constraint (3c) is nonnegative. In model (\ref{WDS}), the quantity $w_{i,t}^\mathrm{CT}$ represents the uncontrollable water load. Similar to the uncontrollable electric load, it is a given value at each period. 

The pumps are considered as constant-speed motors in this paper. The hydraulic characteristics of a constant-speed pump is generally approximated by a quadratic function of the water flow across the pump, i.e. $y^G=a_2f^2+a_1f+a_0$ \cite{UlanickiModeling} and \cite{FooladivandaOptimal}. Contribution of the nonlinear $a_2 f^2$  is usually very small compared to the linear ones $a_1 f+a_0$. Thus, equation \eqref{eq:yf} captures the head gain of a pump simply making $a_2=0$. The following constraints act as the mathematical link between the distribution network (1)--(2) and the WDS (3):
\begin{equation} \label{Pump}
\eta P_{i,t}^{\mathrm{Pump}}=f_{k,t}y_{k,t}^\mathrm{G}=a_{1,k}f_{k,t}^2+a_{0,k}f_{k,t},
\end{equation}
where $i \in \mathcal{N}_\mathrm{E}^\mathrm{P}$ and $k \in \mathcal{E}_\mathrm{W}^\mathrm{P}$.

\subsection{A Co-optimization Framework of Water and Electricity}
Based on the mathematical model of the micro-WEN introduced above, this subsection introduces a co-optimization framework for water and energy networks. The objective of this co-optimization problem is to minimize the total energy cost for meeting the demands of both electricity and water. We formulate the energy cost as
\begin{equation} \label{Objective}
C(P_{i,t}^\mathrm{G})=\sum_t(c_tP_{1,t}^\mathrm{G} + \sum_{i \in \mathcal{N}_E^G/\text{PCC}}(c_{1,i}P_{i,t}^G+c_{2,i}(P_{i,t}^G)^2)),
\end{equation}
where $P_{1,t}^G$ denotes the power from the grid via PCC (i.e. the serial number of PCC is 1); $c_t$ can be considered as the nodal prices at PCC which are obtained by solving the security constrained economic dispatch (SCED) by ISOs/RTOs. As a result, the co-optimization model is
\begin{align}
\begin{split}
\min_{P_{i,t}^G} \quad &(\ref{Objective}) \\
\text{s.t.} \quad &(\ref{PF})-(\ref{Pump})
\end{split}, \tag{CO-OPT}
\end{align}
which is a mixed-integer nonlinear programming (MINLP).

\section{Quasi-Convex Hull Relaxations}

The MINLP problem is computationally intractable, especially for large-scale systems. To reduce the computational burden, this section relaxes the MINLP into a mixed-integer convex programming problem \cite{BonamiAlgorithms} with high-tightness.

\subsection{Convex Hull Relaxations of Constraints (1d) and (2a)}
Within the circular bounds (1e) and (2b), the feasible sets of equations (1d) and (2a) can be captured by the following general formulation:
\begin{equation}
    \Omega_0 = \left\{x \left| \begin{array}{lr}
ax_1^2+bx_2^2=x_3x_4 \\ 
x_1^2+x_2^2 \le c \\
\underline{x}_3,\underline{x}_4 \le x_3,x_4 \le \overline{x}_3,\overline{x}_4
\end{array}\right. \right\} \nonumber
\end{equation}
where $x=[x_1\ x_2\ x_3\ x_4]^T$, $a \geq b$, and $\underline{x}_3\underline{x}_4 \le c \le \overline{x}_3\overline{x}_4$. By generalizing the theorem presented in \cite{LiConvex}, we have the following Lemma.

\textbf{Lemma}. \textit{The convex hull of nonconvex set $\Omega_0$ is}
\begin{equation}
    \Omega_1=CH(\Omega_0) = \left\{x \left| \begin{array}{lr}
ax_1^2+bx_2^2 \le x_3x_4 \\ 
(a-b)x_2^2 + \underline{x}_4x_3 \le ac\\
(a-b)x_2^2 + \underline{x}_3x_4 \le ac\\
D^Tx - d \le 0\\
x_1^2+x_2^2 \le c \\
\underline{x}_3,\underline{x}_4 \le x_3,x_4 \le \overline{x}_3,\overline{x}_4
\end{array}\right. \right\} \nonumber
\end{equation}
\textit{where D = $[0\ 0\ k_1\ k_2]^T$ is a coefficient vector, d is a scalar, and their values are given by}
   \[   
     \begin{cases}         k_1=ac,\,k_2=\underline{x}_3\overline{x}_3,\,d=ac(\underline{x}_3+\overline{x}_3)\ \text{if}\,\overline{x}_3\underline{x}_4\le ac \le \underline{x}_3\overline{x}_4 \\
 k_1=\underline{x}_4\overline{x}_4,\,k_2=ac,\,d=ac(\underline{x}_4+\overline{x}_4)\ \text{if}\,\underline{x}_3\overline{x}_4\le ac \le \overline{x}_3\underline{x}_4 \\      k_1=\underline{x}_4,\,k_2=\underline{x}_3,\,d=ac+\underline{x}_3\underline{x}_4\ \ \text{if}\ \,\overline{x}_3\underline{x}_4,\,\underline{x}_3\overline{x}_4 \le ac  \\
 k_1=\overline{x}_4,\,k_2=\overline{x}_3,\,d=ac+\overline{x}_3\overline{x}_4\ \ \text{if}\ \,ac \le \overline{x}_3\underline{x}_4,\,\underline{x}_3\overline{x}_4
 \end{cases} \nonumber
  \]

\textbf{Proof}: See appendix. $\square$

For the case of (1d), $a=b=1$, $c=\overline{S}_{ik}^2$, and $(\underline{x}_3,\underline{x}_4,\overline{x}_3,\overline{x}_4)=(\underline{\mathcal{V}}_i,0,\overline{\mathcal{V}}_i,\overline{\mathcal{I}}_{ik})$, assume that $\overline{S}_{ik}^2 \le \underline{\mathcal{V}}_i\overline{\mathcal{I}}_{ik}$. Within the system bounds (1e) - (1g), the convex hull of (1d) is given as 
\begin{align} \label{CH1}
\begin{cases}        
     P_{ik,t}^2+Q_{ik,t}^2 \le \mathcal{V}_{i,t}\mathcal{I}_{ik,t} \\
\overline{S}_{ik}^2\mathcal{V}_i+\underline{\mathcal{V}}_i\overline{\mathcal{V}}_i\mathcal{I}_{ik} \le \overline{S}_{ik}^2(\underline{\mathcal{V}}_i+\overline{\mathcal{V}}_i) 
\end{cases}.
\end{align}

For the case of (2a), $a=r_i^{\mathrm{Batt}} + r_i^{Cvt}$, $b=r_i^{Cvt}$, $c=(\overline{S}_{i}^{ES})^2$, and $(\underline{x}_3,\underline{x}_4,\overline{x}_3,\overline{x}_4)=(\underline{\mathcal{V}}_i,0,\overline{\mathcal{V}}_i,+\infty)$.  It is obvious that $\overline{\mathcal{V}}_i\underline{L}_{i,t}^{ES} \le (\overline{S}_{i}^{ES})^2(r_i^{\mathrm{Batt}} + r_i^{Cvt}) \le \underline{\mathcal{V}}_i\overline{L}_{i,t}^{ES}$. Hence, within the system bounds (1g) and (2b), the convex hull of (2a) is
\begin{align}\label{CH2}
\begin{cases}        
     (r_i^{\mathrm{Batt}} + r_i^{Cvt})(P_{i,t}^{ES})^2 + r_i^{Cvt}(Q_{i,t}^{ES})^2 \le L_{i,t}^{ES}\mathcal{V}_{i,t} \\
     r_i^{\mathrm{Batt}}(Q_{i,t}^{ES})^2+\underline{\mathcal{V}}_iL_{i,t}^{ES} \le (\overline{S}_{i}^{ES})^2(r_i^{\mathrm{Batt}} + r_i^{Cvt})\\
(\overline{S}_{i}^{ES})^2\mathcal{V}_i+\underline{\mathcal{V}}_i\overline{\mathcal{V}}_iL_{i,t}^{ES} \le (\overline{S}_{i}^{ES})^2(\underline{\mathcal{V}}_i+\overline{\mathcal{V}}_i) 
\end{cases}.
\end{align}

\subsection{Quasi-Convex Hull Relaxation of (3b)}

\begin{figure}[tb]
\centering
\includegraphics[width=0.33\textwidth]{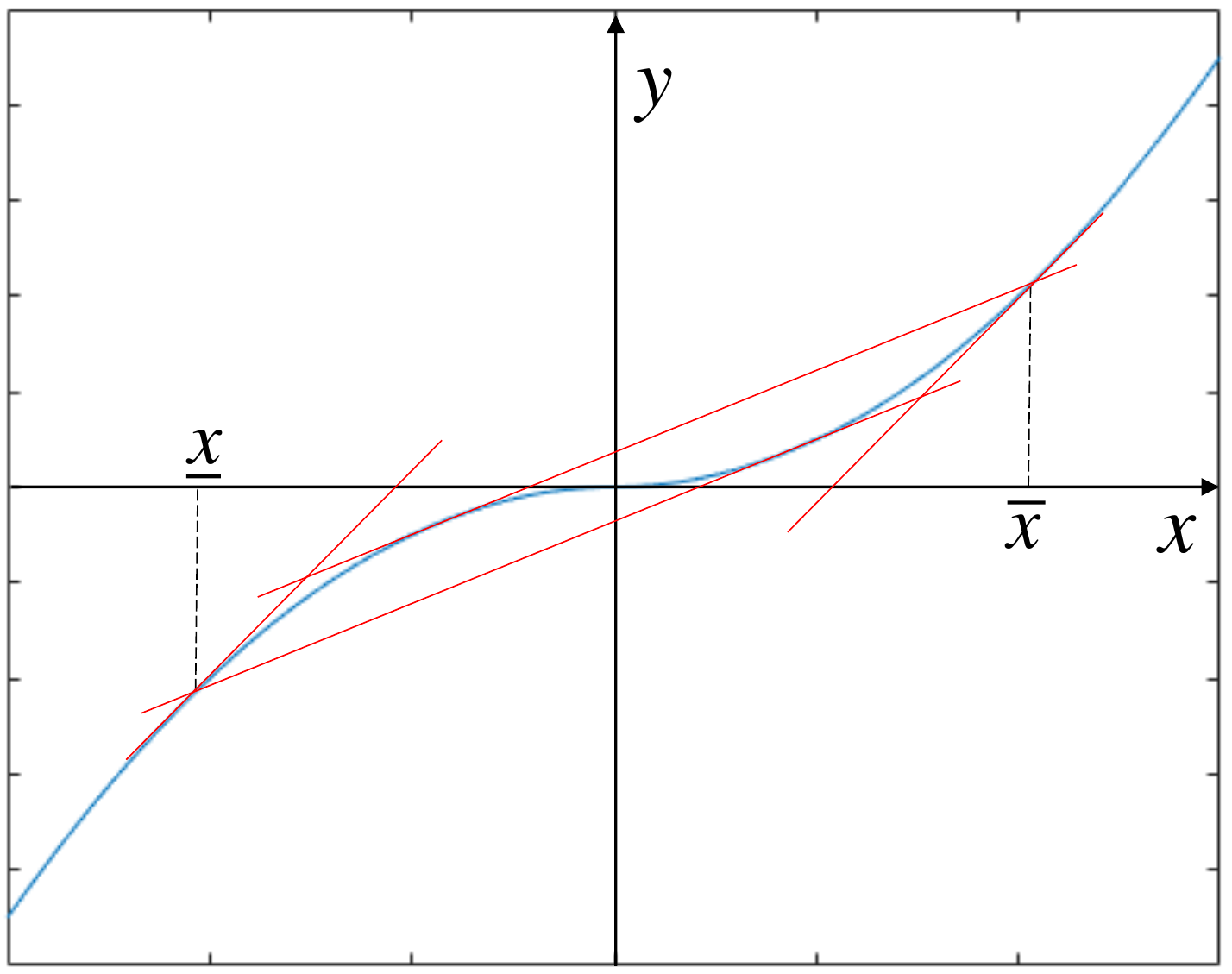}
  \caption{The quasi-convex hull of the hydraulic characteristic of a normal pipe.}
  \label{fig:hydraulic}
\end{figure}

With the left-hand-side replaced by an auxiliary variable $y$ and the right-hand-side replaced by a general term $R^P\mathrm{sgn}(f)f^2$, function (3b) yields the blue curve, as shown in Figure \ref{fig:hydraulic}, in the ($f,y$)-plane. It is relaxed into the red polygon as shown in the figure. The red polygon is not exactly the convex hull of (3b). However, it is very close to the convex hull from the perspective of tightness and, therefore, is called a quasi-convex hull. Its mathematical formulation is given as
\begin{align} \label{PumpCH}
&y_{i,t}-y_{j,t}+h_i-h_j \nonumber \\ 
   &\begin{cases}
\leqslant (2\sqrt{2}-2)R_{k}^P\overline{f}_{k}f_{k,t} +(3-2\sqrt{2})R_{k}^P\overline{f}_{k}^2 \\
\geqslant (2\sqrt{2}-2)R_{k}^P\underline{f}_{k}f_{k,t} +(3-2\sqrt{2})R_{k}^P\underline{f}_{k}^2 \\
\geqslant 2R_{k}^P\overline{f}_{k}f_{k,t} -R_{k}^P\overline{f}_{k}^2 \\
\leqslant 2R_{k}^P\underline{f}_{k}f_{k,t} -R_{k}^P\underline{f}_{k}^2
\end{cases} .
\end{align} 

\subsection{Convex Hull Relaxation of Constraint (3c)}
\begin{figure}[h!]
\centering
\includegraphics[width=0.35\textwidth]{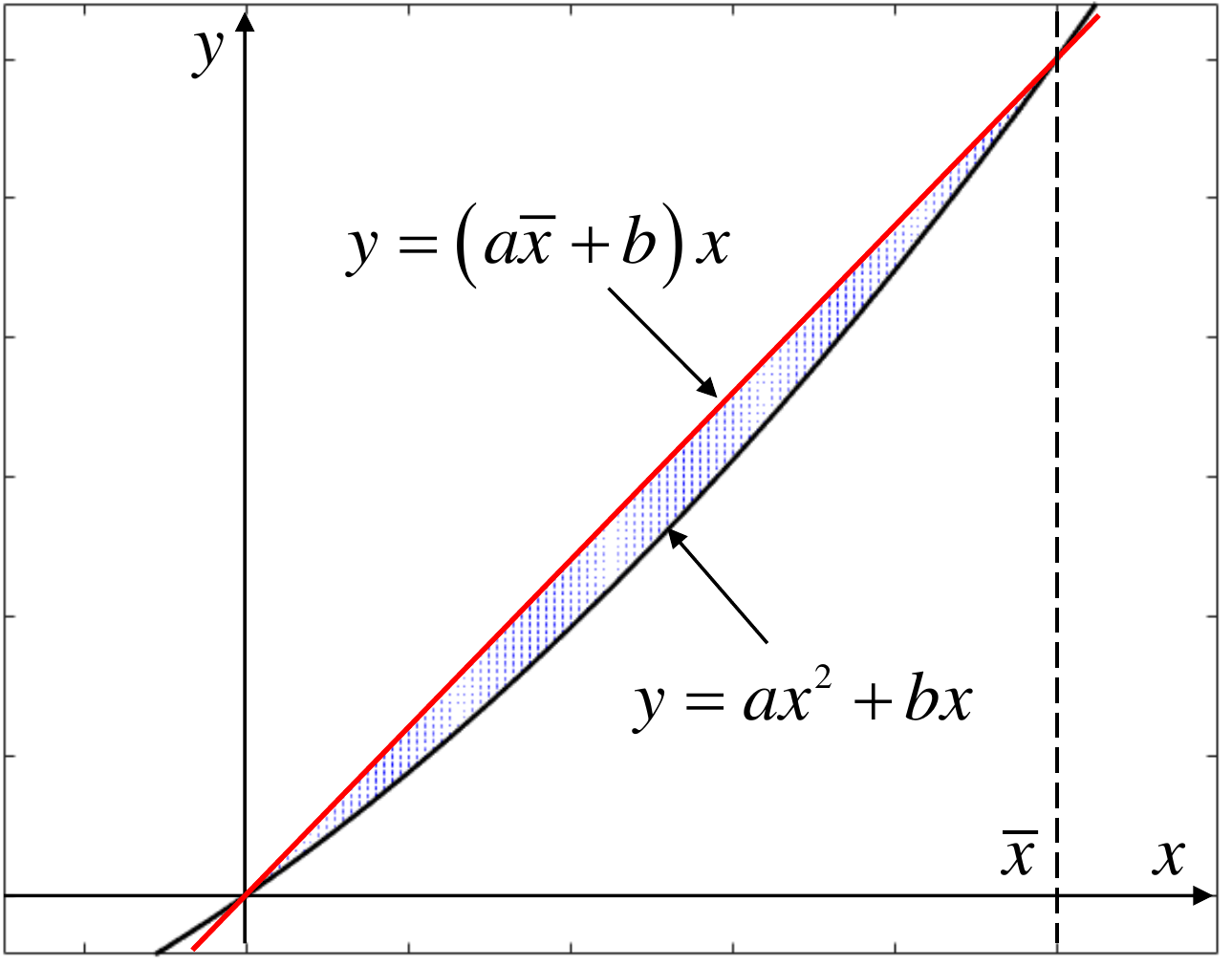}
  \caption{Convex hull of a parabola.}
  \label{fig:parabola}
\end{figure}

The nonconvex constraint (3c) contains a logic proposition. To eliminate the $if$ expression, we use the big-$M$ technique to rewrite constraint (3c) as 
\begin{subequations}
\begin{gather}
R_{k}^Pf_{k,t}^2 - Y_1 \le 0 \label{pumpconvex}\\
Y_2 - R_{k}^Pf_{k,t}^2 \le 0 \label{pumpnonconvex} \\
0 \le f_{k,t} \le M*\alpha 
\end{gather}
\end{subequations}
where $Y_1=y_{i,t}-y_{j,t}+h_i-h_j+y_{k,t}^G+M*(1-\alpha)$ and $Y_2=y_{i,t}-y_{j,t}+h_i-h_j+y_{k,t}^G+M*(\alpha-1)$. Note that the expression (\ref{pumpconvex}) is convex, while (\ref{pumpnonconvex}) is a concave constraint. The convex hull of (\ref{pumpnonconvex}) can be obtained through a geometric approach as shown in Figure \ref{fig:parabola}. Its mathematical expression  is given as 
\begin{equation}
Y_2 - R_{k}^P\overline{f}_{k,t}f_{k,t} \le 0.
\end{equation}

It can be observed, by comparing Figures \ref{fig:hydraulic} and \ref{fig:parabola}, that one can construct a tighter relaxation for the hydraulic characteristic of a pipe if the direction of water flow is given. A planning problem of gas networks is discussed in \cite{BorrazConvex} where the authors introduced additional binary variables and bilinear equations to relax the Weymouth equation. The Weymouth equation is similar to constraint (3b). However, the convex relaxation of the Weymouth equation developed in \cite{BorrazConvex} is not necessarily tighter than the proposed relaxation (\ref{PumpCH}) due to the introduced bilinear equations. Moreover, the auxiliary binary variables and constraints in \cite{BorrazConvex} are not desirable for an operation problem which is sensitive to computational time.

\subsection{Convex Hull Relaxation of Constraint (4)}
Constraint (4) is a quadratic equation which can be considered as the intersection of a convex inequality and a concave inequality. The geometric approach introduced in Figure \ref{fig:parabola} can be used to construct the convex hull of the concave inequality. As a result, the convex hull constraint (4) is given by
\begin{subequations}
\begin{gather}
 \eta P_{i,t}^{\mathrm{Pump}} \ge a_{k}f_{k,t}^2+b_{k}f_{k,t}\\
\eta P_{i,t}^{\mathrm{Pump}} \le (a_{k}\overline{f}_{k,t}+b_{k})f_{k,t}
\end{gather}
\end{subequations}
where $f_{k,t}$ is nonnegative since the direction of pump flows is determined.

\subsection{Quasi-Convex Hull Relaxation of (CO-OPT)}
To sum up, the quasi-convex hull relaxation of the overall co-optimization problem (CO-OPT) is 
\begin{align}
\begin{split}
\min \quad &(\ref{Objective}) \\
\text{s.t.} \quad &\text{(1a-c, e-i), (2b-c), (3a, d-i), }(\ref{CH1}) \\
&(\ref{CH2}),\,(\ref{PumpCH})\text{, (9a, c), (10) and (11)},
\end{split} \tag{C-CO-OPT}
\end{align}
which is a mixed-integer convex quadratically-constrained quadratic programming (MICQCQP) problem.

The basic idea of the quasi-convex hull relaxation of an optimization problem is replacing the nonconvex constraints with their convex hulls or quasi-convex hulls. As discussed in \cite{LiThe,LiConvex}, the concept of convex hull is attractive since the extreme points of a convex hull generally belong to its original non-convex set. If the objective function is a convex function and monotonic over the convex hull, the optimal solution is usually located at one of the extreme points, implying that the optimal solution obtained by solving the convex relaxation is most likely the exact globally optimal solution of the original problem. Unfortunately, for many of the nonlinear nonconvex sets, it is extremely hard to formulate their convex hulls. For such cases, an interesting alternative is to construct a convex inner approximation of a nonconvex set \cite{HungInner}. Compared with the convex relaxation, the foremost advantage of the convex inner approximation is guaranteeing the feasibility of the obtained solutions to the original noncovex set.

A characteristic of the MINLP problem (CO-OPT) is that the integer variables only exist in linear terms of constraints. For the purpose of improving the computational efficiency, it is wise to relax such a mixed-integer problem into a mixed-integer convex problem. The convex relaxations which are tight for the continuous cases are equivalently tight for the discrete cases since the nonconvex terms that need to be relaxed do not contain integer variables.

\section{A Flexible Irrigation Scheme}
\subsection{Flexibility of Irrigation Systems}
It is straightforward to improve the grid flexibility by allowing for controllability of electric loads. From a different angle, this subsection explores opportunities for improving the DR capacity of water systems by investigating the flexibility of customer-owned water tanks and irrigation systems which are water loads rather than electric loads. An intuitive interpretation is that customers can use the superfluous energy to pump and store water. However, there are some certain volume limits on tanks. Thus, it is necessary to develop a coordination strategy for charging and discharging tanks based on the multi-period energy imbalance and the flexible water consumption.

Assuming that the crops growth is not sensitive to the watering time, the irrigation process is considered flexible. To develop a mathematical model of such a flexible irrigation system, we have the following assumptions: i) the irrigation flow is fixed and the irrigation volume is a function of the watering time; ii) for a given season, the total amount of water is fixed, which can be represented as turning on the irrigation system for totally $N$ hours per day. Consequently, the mathematical model is given by
\begin{subequations} \label{irrigation}
\begin{gather}
 \underline{S}_{i}^{w} \leqslant S_{i,0}^w + \sum^t_{t=0} (f^{CT}_{i,t}-f^D_{i,t}-k\alpha_{i,t})  \leqslant \bar{S}_{i}^{w}\\
 \sum_t\alpha_{i,t}=N,
\end{gather}
\end{subequations}
which are mixed-integer linear.

\subsection{A Mixed-integer Convex DSM Scheme of Water Systems}
By incorporating the model (\ref{irrigation}) of flexible irrigation systems into (C-CO-OPT), we have the following DSM scheme
\begin{align}
\begin{split}
\min \quad &(\ref{Objective}) \\
\text{s.t.} \quad &\text{(1a-c, e-i), (2b-c), (3a, d-i), }(\ref{CH1}) \\
&(\ref{CH2}),\,(\ref{PumpCH})\text{, (9a, c), (10), (11) and (\ref{irrigation})},
\end{split} \tag{C-DSM}
\end{align}
which is also an MICQCQP problem.

\section{Case Study}

\subsection{Introduction to the Test System}
\begin{figure}[t]
\centering
\includegraphics[width=0.5\textwidth]{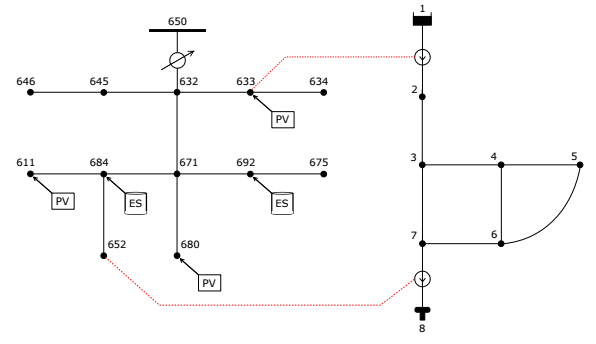}
  \caption{Topology of the test system.}
  \label{fig:Test}
\end{figure}
\begin{figure}[t]
\centering
\includegraphics[width=0.5\textwidth]{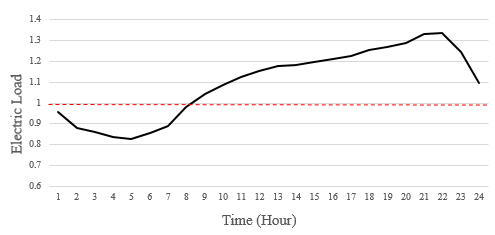}
  \caption{The shape of a typical average load in summer.}
  \label{fig:load}
\end{figure}
The micro-WEN for the case study is composed of the IEEE 13-bus system and an 8-node WDS from the EPANET manual \cite{EPANET}. The topology of the test micro-WEN is given in Figure \ref{fig:Test}. We assume that the 13-bus microgrid is integrated with high penetration of PV resources. Figure 4 shows the shape of a typical average load in summer \cite{NarangHigh}. The 24-hour load profile of the 13-bus system is generated by applying this load shape with the load provided by IEEE as the load at 9 am. Further detailed information about the PV systems, BESSs, and pumps are given in Table \ref{Table1}.
\begin{table}[h]
\centering
\caption{PV system, BESS, and Pump Parameters}
\label{Table1}
\begin{tabular}{cc}
\hline\hline
\textbf{PV location (bus \#) and capacity} & \textbf{Penetration} \\
633 (0.5 MW), 680 (0.2 MW), 684 (0.5 MW)   & 34.68\%               \\
\hline
\multicolumn{2}{c}{\textbf{BESS location (bus \#) and capacity}}   \\
\multicolumn{2}{c}{684 (1.15 MVA, 2.5 MWh), 692 (1.41 MVA, 3.2 MWh)} \\
\hline
\multicolumn{2}{c}{\textbf{Pump location (bus \#) and parameters}}   \\
\multicolumn{2}{c}{633 (\textit{b}=0.3 p.u., \textit{c}=0.4 p.u.), 652 (\textit{b}=0.3 p.u., \textit{c}=0.4 p.u.)} \\
\hline\hline
\end{tabular}
\end{table}

Pumps deliver 30.48 meters and 15.24 meters of head respectively at a flow of 0.038 m$^3$/s. The tank is 18.3 meters in diameter and 5.1 meters in depth. For the 24-hour demand profiles and the lengths of pipes of the water system, please refer to the EPANET manual. The parameter $R_{ij}^P$ is calculated by (the subscript $ij$ is eliminated for the sake of simplicity)
\begin{align}
R^P=\frac{8fL}{\pi^2gD^5} \nonumber
\end{align}
where $f$ is the coefficient of surface resistance, $D$ and $L$ are the diameter and length of pipe respectively, and $g$ is the gravitational acceleration.

\subsection{Tightness of the QCH Relaxation}
The tightness of the proposed convex relaxation is first evaluated by comparing solutions obtained by solving (CO-OPT) and (C-CO-OPT), respectively. Using the JuMP package of Julia \cite{KwonJulia}, the optimization problems were solved in a MAC computer with a 64-bit Intel i7 dual core CPU at 2.40 GHz and 8 GB of RAM. The MINLP problem (CO-OPT) and its quasi-convex hull relaxation (C-CO-OPT) are solved by calling BONMIN \cite{Grossmann} and GUROBI \cite{GUROBI} solvers respectively. 

The simulation results are tabulated in Table \ref{tightness}. The first and foremost improvement brought by convexification is in the computational efficiency. The required CPU time has been significantly reduced by solving the quasi-convex hull relaxation. Note that BONMIN is an open source solver with a limited computational capacity. However, MINLP problems are NP-hard to solve. Even using some well-designed commercial solvers, like KNITRO \cite{knitro}, the CPU time of solving a MINLP problem is still not comparable to that of solving a MICQCQP problem of a similar size.

The proposed quasi-convex hull relaxation is exact for the test case in this paper. It means that the optimal solution obtained by solving (C-CO-OPT) is the exact global optimal solution of its original nonconvex problem (CO-OPT) with zero optimality gap. The numerical results in this subsection demonstrate the potential of the proposed quasi-convex hull relaxation for convexifying similar MINLP problems with high-accuracy. 

\begin{table*}[ht]
\centering
\caption{Results of Tightness}
\label{tightness}
\begin{tabular}{cccccc}
\hline\hline
\textbf{Problem} & \textbf{\begin{tabular}[c]{@{}c@{}}Mathematical\\ Classification\end{tabular}} & \textbf{Solver} & \textbf{\begin{tabular}[c]{@{}c@{}}Optimal\\ Solution (\$)\end{tabular}} & \textbf{\begin{tabular}[c]{@{}c@{}}Optimality\\ Gap\end{tabular}} & \textbf{CPU Time} \\
\hline
(CO-OPT)         & MINLP                                                                          & BONMIN          & 1463                                                                     & -                                                                 & $\approx$ 2 hrs             \\
(C-CO-OPT)       & MICQCQP                                                                        & GUROBI          & 1463                                                                     & 0                                                                 & $<$ 1 sec  \\
\hline\hline
\end{tabular}
\end{table*}

\begin{figure}[t]
\centering
\includegraphics[width=0.5\textwidth]{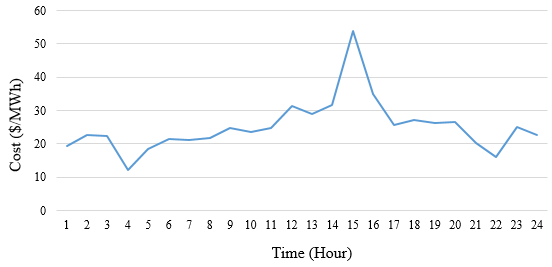}
  \caption{The 24-hour nodal price at PCC.}
  \label{fig:LMP}
\end{figure}

\subsection{Improvement in System Security}
In current practice, the electrical and water systems are operated separately by the electrical and water utilities respectively. First, the water utility tries to minimize the energy consumption by doing a day-ahead optimal pump scheduling based on the day-ahead water demand forecast (see the formulation (OPS) in \cite{LiModeling}). Then, with the schedule of energy consumption reported by the water utility, the power operator solves a multi-period optimization problem (see the formulation (UC) in \cite{LiModeling}) to minimize the total energy cost for meeting all electricity demands based on the day-ahead forecast of electricity demands, where the diesel generators and BESS units are control devices. It can be observed from Table \ref{Security} that the co-optimization produces a lower-cost solution under the same penetration level.

In this section, we also evaluate the improvement in system security introduced by co-optimizing the electrical and water networks. There is a certain limit on the penetration of PV that the distribution system can accommodate. A PV generation that exceed this limit will cause security problem to the system. The simulation results tabulated in Table \ref{Security} demonstrate that, by considering the pumps as controllable loads, the PV penetration which can be accommodated by the power distribution system is increased by 33$\%$. In reality, the proposed co-operation scheme of micro-WEN allows the  electricity-driven water facilities response to energy uncertainties in the power grid and, consequently, improve the system's security.
\begin{table}[h!]
\centering
\caption{Capability of PV Penetration}
\label{Security}
\begin{tabular}{ccc}
\hline\hline
\textbf{Operation Scheme}   & \textbf{Penetration Cap} & \textbf{Optimal Solution (\$)} \\
\hline
Independent Optimization & 83.23\%  & 1579                    \\
Co-optimization    & 110.98\%  & 1463 \\
\hline\hline
\end{tabular}
\end{table}

\subsection{Efficiency of the DSM Scheme of Tanks}
This subsection evaluates the efficiency for the demand response of the flexible irrigation system (\ref{irrigation}) by comparing the results of (C-CO-OPT) and (C-DSM). Both problems are MICQCQP and solved by GUROBI on the computer mentioned in the previous subsection. The results are tabulated in Table \ref{Economic}. Problem (C-DSM) has more integer decision variables than (C-CO-OPT). However, the CPU time for solving (C-DSM) is not significantly larger than that for solving (C-CO-OPT) due to the property that the nonlinear terms are convex. 

In the studied case, 30$\%$ of the total water load is for irrigation and, namely, flexible. It can be observed from Table \ref{Economic} that the total operational costs of the mirco-WEN can be reduced by considering the flexibility of the irrigation systems. A considerable cost saving can be expected if: 1) the proposed approach is applied to a larger system, 2) the flexibility of other water facilities, such as water/waste water treatment, desalination, recycling, and cooling, is also included.

\section{Conclusion and Future Work}
This paper introduces a mixed-integer nonlinear mathematical model for the distribution-level water-energy nexus, or the micro-WEN. A co-optimization framework for water and energy distribution networks is built upon this mathematical model. Based on the convex hulls or quasi-convex hulls of the system components, a tight mixed-integer convex relaxation is developed to improve the computational efficiency of solving the co-optimization framework. When the proposed approach is applied to solve the co-optimization problem of a micro-WEN which consists of a 13-bus distribution system and an 8-node water distribution network, the CPU time reduces from nearly 2 hours to less than 1 second. 

To further explore the capacity of water distribution systems for providing demand response service to the power grids, this paper developed an optimal demand response framework considering a flexible irrigation system. Simulation results on the test micro-WEN demonstrated the DR potential of water systems. In the future work, we will consider the flexibility of other water facilities, such as water/waste water treatment, desalination, recycling, and cooling, in the DSM model.

The water-energy nexus is a fundamental infrastructure in the building/city/village as both water and electricity are lifelines of humans. The findings of this research should be a good fit to the research framework of the smart building/city/village.

\begin{table}[h]
\centering
\caption{Results of Economic Efficiency}
\label{Economic}
\begin{tabular}{ccccc}
\hline\hline
\textbf{Problem} & \textbf{\begin{tabular}[c]{@{}c@{}}Mathematical\\ Classification\end{tabular}} & \textbf{Solver} & \textbf{\begin{tabular}[c]{@{}c@{}}Optimal\\ Solution (\$)\end{tabular}} & \textbf{CPU Time} \\
\hline
(C-CO-OPT)       & MICQCQP                                                                        & GUROBI          & 1463                                                                     & $<$ 1 sec            \\
(C-DSM)          & MICQCQP                                                                        & GUROBI          & 1455                                                                     & $<$ 1 sec  \\
\hline\hline
\end{tabular}
\end{table}


%

\appendix
The proof for only the case, where $k_1=ac,\,k_2=\underline{x}_3\overline{x}_3$, and $d=ac(\underline{x}_3+\overline{x}_3)$ is provided due to the page limit. The set $\Omega_1$ represents a convex solid, in the $x$-space, that consists of 5 (linear) facets and 4 (nonlinear) surfaces. The relation $\Omega_1=CH(\Omega_0)$ means $CH(\Omega_0)\subseteq \Omega_1$ and $\Omega_1 \subseteq CH(\Omega_0)$. \\
\textbf{(i) $CH(\Omega_0)\subseteq \Omega_1$} 

\begin{figure}[h!]
\centering
\includegraphics[width=0.35\textwidth]{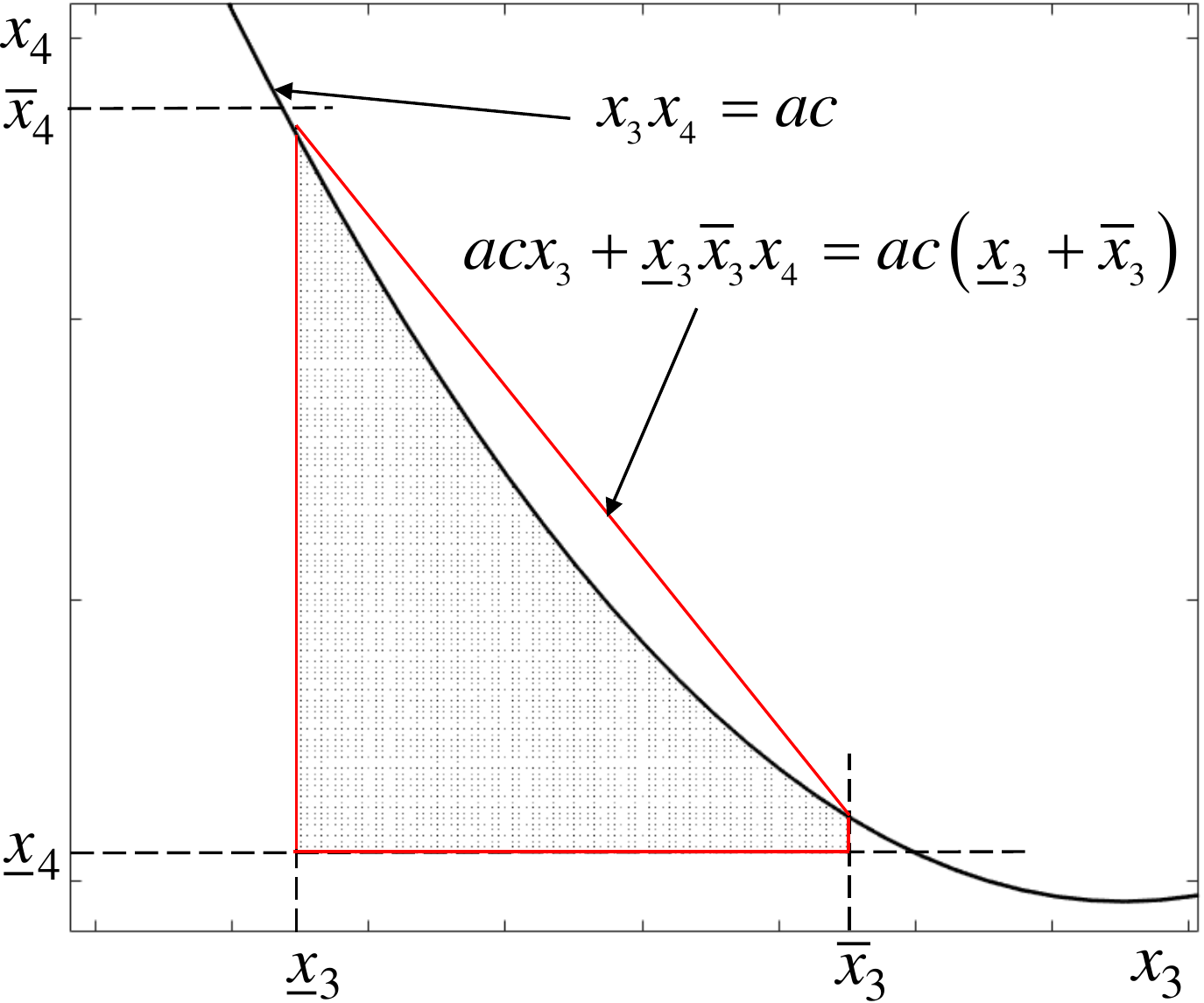}
  \caption{Convex hull relaxation of $x_3x_4 \le ac$. The shaded area is the original nonconvex set while the trapezoid region with red boundaries is its convex hull.}
  \label{fig:proof}
\end{figure}

Let
\begin{equation}
    \Omega_2= \left\{x \left| \begin{array}{lr}
ax_1^2+bx_2^2 \le x_3x_4 \\ 
x_1^2+x_2^2 \le c \\
\underline{x}_3,\underline{x}_4 \le x_3,x_4 \le \overline{x}_3,\overline{x}_4
\end{array}\right. \right\} \nonumber
\end{equation}
\begin{equation}
    \Omega_3= \left\{x \left| \begin{array}{lr}
 (a-b)x_2^2 + \underline{x}_4x_3 \le ac\\
x_1^2+x_2^2 \le c \\
\underline{x}_3,\underline{x}_4 \le x_3,x_4 \le \overline{x}_3,\overline{x}_4
\end{array}\right. \right\} \nonumber
\end{equation}
\begin{equation}
    \Omega_4 = \left\{x \left| \begin{array}{lr}
(a-b)x_2^2 + \underline{x}_3x_4 \le ac\\
x_1^2+x_2^2 \le c \\
\underline{x}_3,\underline{x}_4 \le x_3,x_4 \le \overline{x}_3,\overline{x}_4
\end{array}\right. \right\} \nonumber
\end{equation}
\begin{equation}
    \Omega_5 = \left\{x \left| \begin{array}{lr}
D^Tx - d \le 0\\
x_1^2+x_2^2 \le c \\
\underline{x}_3,\underline{x}_4 \le x_3,x_4 \le \overline{x}_3,\overline{x}_4
\end{array}\right. \right\}, \nonumber
\end{equation}
which are convex sets. It is straight forward to know that $\Omega_0 \subseteq \Omega_2$. Under the condition $x_1^2+x_2^2 \le c$, equation $ax_1^2+bx_2^2 = x_3x_4$ implies $(a-b)x_2^2 + x_3x_4 = a(x_1^2+x_2^2) \le ac$. Any point $(x_2,x_3)$ that satisfies $(a-b)x_2^2 + x_3x_4 \le ac$ will also satisfy $(a-b)x_2^2 + \underline{x}_4x_3 \le ac$. Therefore, $\Omega_0 \subseteq \Omega_3$. Similarly, we can prove that $\Omega_0 \subseteq \Omega_4$. Equation $ax_1^2+bx_2^2 = x_3x_4$ also implies $x_3x_4 \le ac$ since $ax_1^2+bx_2^2 \le ac$. From Figure \ref{fig:proof}, it suffices to show that $\Omega_0 \subseteq \Omega_5$. The convex hull of $\Omega_0$ is defined as the intersection of all convex relaxations of $\Omega_0$ \cite{Boyd}. Hence, $CH(\Omega_0) \subseteq (\Omega_2 \cap \Omega_3 \cap \Omega_4 \cap \Omega_5) = \Omega_1$.

\begin{flushleft}
\textbf{(ii) $\Omega_1 \subseteq CH(\Omega_0)$} 
\end{flushleft}

If a linear cut is valid for the convex set $CH(\Omega_0)$, it will also be valid for any subset of $CH(\Omega_0)$. Note that “a linear inequality is valid for a set” means the inequality is satisfied by all its feasible solutions \cite{Borozan}. On the other hand, $\Omega_1$ is a subset of $CH(\Omega_0)$ if any valid cut of $CH(\Omega_0)$ is also valid for $\Omega_1$ according to the properties of supporting hyperplanes \cite{Boyd}. Let $\alpha^Tx \le \beta$ denote any given valid linear cut for $CH(\Omega_0)$, it should also be valid for $\Omega_0$. The cut $\alpha^Tx \le \beta$ is valid for $\Omega_1$ if it is valid for all the surfaces of $\Omega_1$.

The mathematical formulation of a surface can be obtained by changing one inequality constraint in $\Omega_1$ into an equality. For instance, $\mathrm{Surf}_1$ and $\mathrm{Surf}_2$ are two surfaces of solid $\Omega_1$. It is straight forward to show that $\alpha^Tx \le \beta$ is valid for $\mathrm{Surf}_1$ since, in reality, $\mathrm{Surf}_1=\Omega_0$. In this appendix, we prove that $\alpha^Tx \le \beta$ is valid for $\mathrm{Surf}_2$ as an example.
\begin{equation}
    \mathrm{Surf}_1 = \left\{x \left| \begin{array}{lr}
ax_1^2+bx_2^2 = x_3x_4 \\ 
(a-b)x_2^2 + \underline{x}_4x_3 \le ac\\
(a-b)x_2^2 + \underline{x}_3x_4 \le ac\\
D^Tx - d \le 0\\
x_1^2+x_2^2 \le c \\
\underline{x}_3,\underline{x}_4 \le x_3,x_4 \le \overline{x}_3,\overline{x}_4
\end{array}\right. \right\} \nonumber
\end{equation}
\begin{equation}
    \mathrm{Surf}_2 = \left\{x \left| \begin{array}{lr}
ax_1^2+bx_2^2 \le x_3x_4 \\ 
(a-b)x_2^2 + \underline{x}_4x_3 = ac\\
(a-b)x_2^2 + \underline{x}_3x_4 \le ac\\
D^Tx - d \le 0\\
x_1^2+x_2^2 \le c \\
\underline{x}_3,\underline{x}_4 \le x_3,x_4 \le \overline{x}_3,\overline{x}_4
\end{array}\right. \right\} \nonumber
\end{equation}

Let $\tilde{x}=(\tilde{x}_1,\tilde{x}_2,\tilde{x}_3,\tilde{x}_4)$ denote any given point on $\mathrm{Surf}_2$, where $(a-b)\tilde{x}_2^2 + \underline{x}_4\tilde{x}_3 = ac$. Let's consider the two chosen points $\hat{x}_1=(\hat{x}_{1,1},\tilde{x}_2,\tilde{x}_3,\hat{x}_{1,4})$ and $\hat{x}_2=(\hat{x}_{2,1},\tilde{x}_2,\tilde{x}_3,\hat{x}_{2,4})$ which are located in the original feasible set $\Omega_0$. That means $a\hat{x}_{i,1}^2+b\tilde{x}_2^2=\tilde{x}_3\hat{x}_{i,4}$ for $i=1,2$. By carefully choosing the values of $\hat{x}_{1,1}$ and $\hat{x}_{2,1}$ inside the required bounds, we can make conditions $\hat{x}_{1,1} \le \tilde{x}_1 \le \hat{x}_{2,1}$ and $\hat{x}_{1,4} \le \tilde{x}_4 \le \hat{x}_{2,4}$ hold. As a result, it suffices to show that the following condition holds:
\begin{align}
\tilde{x}=\gamma_1\hat{x}_1+\gamma_2\hat{x}_2 \nonumber
\end{align}
where $\gamma_1$ and $\gamma_2$ are nonnegative, and $\gamma_1+\gamma_2=1$. There always exist two points $\hat{x}_1$ and $\hat{x}_2$ that satisfy the above conditions for any given point $\tilde{x}$ on $\mathrm{Surf}_2$ as long as the intersection of $\Omega_0$ and $\mathrm{Surf}_2$ is nonempty. In reality, $\mathrm{Surf}_2$ is redundant if its intersection with $\Omega_0$ is an empty set. We don't need to consider the case that $\mathrm{Surf}_2$ is redundant.

Given that $\hat{x}_1$ and $\hat{x}_2$ belong to the original set $\Omega_0$, the linear cut $\alpha^Tx \le \beta$ is valid for both of them. Therefore, we have
\begin{align}
\alpha^T\tilde{x}=\gamma_1\alpha^T\hat{x}_1+\gamma_2\alpha^T\hat{x}_2 \le \gamma_1\beta + \gamma_2\beta=\beta, \nonumber
\end{align}
which means $\alpha^Tx \le \beta$ is also valid for $\tilde{x}$ and, consequently, $\mathrm{Surf}_2$ since $\tilde{x}$ represents any point on $\mathrm{Surf}_2$. We do not provide the proof showing that the linear cut is also valid for the rest of surfaces due to the page limit. The readers could to the proofs themselves by following the above method for the rest of surfaces of $\Omega_1$ and the rest of the cases in the Lemma.
\ifCLASSOPTIONcaptionsoff
  \newpage
\fi

\end{document}